\documentstyle[amscd,amssymb,verbatim,diagrams,12pt]{amsart}
\newcommand{\ZZ}{{\cal Z}}

\renewcommand{\mod}{\operatorname{mod}}

\newcommand{\OO}{{\cal O}}

\newcommand{\coker}{\operatorname{coker}}

\newcommand{\be}{{\bf e}}

\newcommand{\G}{{\Bbb G}}

\newcommand{\mg}{{\frak m}}
\newcommand{\hra}{\hookrightarrow}
\newcommand{\lan}{\langle}
\newcommand{\ran}{\rangle}

\newcommand{\UU}{{\cal U}}
\newcommand{\WW}{{\cal W}}

\newcommand{\Spec}{\operatorname{Spec}}

\newcommand{\Proj}{\operatorname{Proj}}
\renewcommand{\P}{{\Bbb P}}

\newcommand{\A}{{\Bbb A}}

\numberwithin{equation}{subsection}

\newtheorem{thm}{Theorem}[subsection]
\newtheorem{prop}[thm]{Proposition}
\newtheorem{lem}[thm]{Lemma}
\newtheorem{cor}[thm]{Corollary}
{  \theoremstyle{definition}
\newtheorem{defi}[thm]{Definition}

\newtheorem{rem}[thm]{Remark}

}

\newcommand{\Pf}{\noindent {\it Proof}}

\newcommand{\ov}{\overline}

\newcommand{\cusp}{\operatorname{cusp}}

\newcommand{\FF}{{\cal F}}

\newcommand{\MM}{{\cal M}}

\newcommand{\VV}{{\cal V}}

\newcommand{\Exc}{\operatorname{Exc}}

\renewcommand{\a}{\alpha}
\renewcommand{\b}{\beta}

\newcommand{\De}{\Delta}
\newcommand{\la}{\lambda}

\newcommand{\C}{{\Bbb C}}

\newcommand{\Z}{{\Bbb Z}}
\newcommand{\Q}{{\Bbb Q}}

\newcommand{\wt}{\widetilde}

\newcommand{\sub}{\subset}
\newcommand{\ed}{\qed\vspace{3mm}}

\newcommand{\forg}{\operatorname{for}}

\newcommand{\ba}{{\bf a}}
\newcommand{\dashra}{\dashrightarrow}

\title{Contracting the Weierstrass locus to a point}
\author{Alexander Polishchuk}
\thanks{Supported in part by the NSF grant DMS-1400390}

\begin{document}
\maketitle

\begin{abstract}
We construct an open substack $U\sub\MM_{g,1}$ with the complement of codimension $\ge 2$
and a morphism from $U$ to a weighted projective stack, which sends the Weierstrass locus $\WW\cap U$ to a point,
and maps $\MM_{g,1}\setminus\WW$ isomorphically to its image. The construction uses alternative birational models of
$\MM_{g,1}$ and $\MM_{g,2}$ from \cite{P-krich}.
\end{abstract}

\section*{Introduction}

Let $\WW\sub \MM_{g,1}$ denote the locus in the moduli stack of smooth one-pointed curves of genus $g$,
consisting of $(C,p)$ such that $p$ is a Weierstrass point on $C$,
i.e., $h^1(gp)\neq 0$. It is well known that $\WW$ is an irreducible divisor.
In this paper we construct a rational map from $\MM_{g,1}$ to a proper DM-stack with projective coarse moduli space, 
which contracts $\WW$ to a single point and maps $\MM_{g,1}\setminus\WW$ isomorphically to its image
(see Theorem A below). This is partly motivated by the question whether the class of the
closure of $W$ in $\ov{M}_{g,1}$ generates an extremal ray (we do not solve this; however, see Prop.\ \ref{bir-contr-prop},
Rem.\ \ref{bir-contr-rem} and the discussion below). Note that for small $g$ some pointed Brill-Noether divisors
were shown to generate extremal rays in the effective cone of $\ov{M}_{g,1}$ in \cite{Rulla}, \cite{Jensen1} and 
\cite{Jensen2}.

The construction involves certain moduli stacks studied in \cite{P-krich}.
Namely, in \cite{P-krich} we introduced and studied 
the moduli stack of curves with marked points $(C,p_1,\ldots,p_n)$, where $C$ is
a reduced projective curve of arithmetic genus $g$,
such that $h^1(a_1p_1+\ldots+a_np_n)=0$ for fixed integer weights $a_i\ge 0$ such that
$a_1+\ldots+a_n=g$ (we assume that the marked points are smooth and distinct). 
We denote this stack by $\UU^{ns}_{g,n}(a_1,\ldots,a_n)$. We showed that $\UU^{ns}_{g,n}(a_1,\ldots,a_n)$
can be realized as a quotient of an affine scheme by a torus action and studied
the related GIT picture which leads to interesting projective birational models of $M_{g,n}$.
In particular, for $n=1$ and $a_1=g$ there is a unique nonempty GIT quotient stack $\ov{\UU}^{ns}_{g,1}(g)$,
obtained from $\UU^{ns}_{g,1}(g)$ by deleting one point corresponding to the most singular cuspidal curve. Furthermore,
$\ov{\UU}^{ns}_{g,1}(g)$ is a closed substack in a weighted projective space
(see Sec.\ \ref{background-sec} for details).

We start by considering the natural rational map
\begin{equation}\label{rat-map-eq}
\forg_2:\UU^{ns}_{g,2}(g-1,1)\dashra \ov{\UU}^{ns}_{g,1}(g)
\end{equation}
given by forgetting the second marked point
(more precisely, the map $\forg_2$ is regular on a certain open substack which is dense in the component corresponding
to smoothable curves). Our main technical result is that \eqref{rat-map-eq} is regular on the open substack
of $(C,p_1,p_2)$ such that $h^1((g+1)p_1)=0$, and that the divisor, defined by the condition $h^1(gp_1)\neq 0$, gets
contracted to a point (see Prop.\ \ref{reg-map-prop}). Furthermore, we show that this point has trivial group of
automorphisms.
We derive from this the following result.


\bigskip

\noindent
{\bf Theorem A}. {\it Assume that $g\ge 2$. The natural open embedding of stacks
$$\MM_{g,1}\setminus\WW\hra  \ov{\UU}^{ns}_{g,1}(g)$$
extends to a regular morphism
$$\phi=\phi_g:U\to \ov{\UU}^{ns}_{g,1}(g),$$
for some open substack $U\sub \MM_{g,1}$ containing $\MM_{g,1}\setminus\WW$
and such that $\MM_{g,1}\setminus U$ has codimension $\ge 2$ in $\MM_{g,1}$. Furthermore, $\phi$
contracts $U\cap\WW$ to a single point, which has no nontrivial automorphisms.}

\medskip


More precisely, the open substack $U$ in the above Theorem consists of $(C,p)$ such that $h^1((g+1)p)=0$ and
$h^0((g-1)p)=1$.

We study the case $g=2$ in more detail. In this case we get a more precise result involving a certain 
modular compactification of $\MM_{2,1}$. 

Recall that Smyth introduced in \cite{Smyth-modular} 
the notion of an extremal assignment, which is a rule associating to each stable curve of given
arithmetic genus some of its irreducible components (this rule should be stable under degenerations).
For each extremal assignment $\ZZ$, Smyth considered the moduli stack $\ov{\MM}_{g,n}(\ZZ)$
of $\ZZ$-stable curves, i.e., pointed curves 
$C$ for which there exists a stable curve $C'$ and a map of pointed curves $C'\to C$, contracting precisely the components
of $C'$, assigned by $\ZZ$, in a certain controlled way. In this paper we consider only one extremal assignment which associates
to every stable curve all of its unmarked components (see \cite[Ex.\ 1.12]{Smyth-modular}), so when we say $\ZZ$-stable we always mean this particular extremal assignment.

We prove that the map $\phi_2$ extends to a regular morphism of stacks
$$\phi_2:\ov{\MM}_{2,1}(\ZZ)\to \ov{\UU}^{ns}_{2,1}(2)$$
contracting the closure of $\WW$ to one point (see Theorem \ref{g2-thm}). 
Furthermore, we identify the point $\phi_2(\WW)$ explicitly
as a certain cuspidal curve $C_0$ (see Definition \ref{C0-defi}), 
and show that $\phi_2$ induces an isomorphism of the complement of
$\WW$ to the complement of $\phi_2(\WW)$.

We also prove that the natural rational map of the coarse moduli spaces $\ov{M}_{2,1}\dashra \ov{U}^{ns}_{2,1}(2)$
is a birational contraction with the exceptional divisors $\ov{W}$ and $\De_1$ (see Proposition \ref{bir-contr-prop}).
One can expect that the rational map $\ov{M}_{g,1}\dashra \ov{U}^{ns}_{g,1}(g)$ is still a birational contraction
for $g>2$ (see Remark \ref{bir-contr-rem} for further discussion).

In addition, in Sec.\ \ref{g-2-n-1-a-2-sec} we obtain an isomorphism
$$\ov{\UU}^{ns}_{2,1}(2)\simeq \P(2,3,4,5,6),$$
where the right-hand side is the weighted projective stack.


\medskip

\noindent
{\it Conventions}. In Sec.\ \ref{g-2-n-1-a-2-sec} we work over $\Z[1/6]$. Everywhere else
we work over $\C$. By a curve we mean a connected reduced projective curve.
By the genus of a curve we always mean arithmetic genus.
For DM-stacks whose notation involves calligraphic letters $\MM$, $\UU$ and $\WW$, we denote their coarse moduli spaces
by replacing these letters by $M$, $U$ and $W$.

\section{Rational maps $\forg_2$ and $\phi$}\label{rat-map-sec}

\subsection{Moduli spaces of curves with non-special divisors}\label{background-sec}

We start by recalling some results from \cite{P-krich} about the stacks $\UU^{ns}_{g,n}(\ba)$, where
$\ba=(a_1,\ldots,a_n)$ and $a_i$ are non-negative integers with $a_1+\ldots+a_n=g$.
We denote by $\wt{\UU}^{ns}_{g,n}(\ba)$ the $\G_m^n$-torsor over $\UU^{ns}_{g,n}(\ba)$,
corresponding to choices of nonzero tangent vectors at the marked points. It is proved in \cite{P-krich} that
$\wt{\UU}^{ns}_{g,n}(\ba)$ is an affine scheme of finite type.
In this paper we only need the case when all $a_i$ are positive, so we assume this is the case.

The key result we will use is that for each $i=1,\ldots,n$, and each $(C,p_1,\ldots,p_n,v_1,\ldots,v_n)$ in
$\wt{\UU}^{ns}_{g,n}(\ba)$ (where $v_i$ is a nonzero tangent vector at $p_i$), there is a canonical formal parameter
$t_i$ on $C$ at $p_i$, such that $\lan v_i,dt_i\ran=1$, which is defined as follows. Given a formal parameter $t_i$,
for each $m>a_i$ there is unique, up to adding a constant, rational function
$f_i[-m]\in H^0(C,\OO(mp_i+\sum_{j\neq i}a_jp_j))$ with the Laurent expansion in $t_i$ of the form
\begin{equation}\label{f-i-m-eq}
f_i[-m]=t_i^{-m}+\sum_{q\ge -a_i}\a_i[-m,q]t_i^q.
\end{equation}
The canonical parameter is uniquely characterized by the condition that $\a_i[-m,-a_i]=0$ for every $m>a_i$.
Using these formal parameters we can consider for every pair $(i,j)$ and $m>a_i$ the expansion of $f_i[-m]$
at $p_j$:
$$f_i[-m]=\sum_{q\ge -a_j}\a_{ij}[-m,q]t_i^q$$
(note that $\a_i[-m,q]=\a_{ii}[-m,q]$).
Now we can view the coefficients $\a_{ij}[-m,q]$ as functions on $\wt{\UU}^{ns}_{g,n}(\ba)$,
where we fix the ambiguity in adding a constant to $f_i[-m]$ by requiring that $\a_i[-m,0]=0$. It follows from the
results of \cite{P-krich} that these functions are all expressed in terms of a finite number of them, which gives a
closed embedding of $\wt{\UU}^{ns}_{g,n}(\ba)$ into an affine space.

The rescaling of the tangent vectors $(v_i)$ defines an action of $\G_m^n$ on $\wt{\UU}^{ns}_{g,n}(\ba)$, so
that the weight of the function $\a_{ij}[-m,q]$ is $m\be_i+q\be_j$, where $(\be_i)$ is the standard basis in
the character lattice of $\G_m^n$.

There is a special point in $\wt{\UU}^{ns}_{g,n}(\ba)$ which is a unique point stable under the action of
$\G_m^n$: it is the point where all the functions $\a_{ij}[-m,q]$ vanish, i.e., it corresponds to the origin in the
ambient affine space. The underlying curve is the union of $n$ rational cuspidal curves $C^{\cusp}(a_i)$, glued transversally
at the cusp. Here $C^{\cusp}(a)$ is the projective curve with the affine part given by
$\Spec(k\cdot 1+x^{a+1}k[x])$, and with one smooth point at infinity (see \cite[Sec.\ 2.1]{P-krich}).

In \cite{P-krich} we also studied the GIT picture for the $\G_m^n$-action on $\wt{\UU}^{ns}_{g,n}(\ba)$. 
In general we have stability conditions depending on a character $\chi$ of $\G_m^n$. 
In the case $n=1$, i.e., for $\wt{\UU}^{ns}_{g,1}(g)$ there is a unique nonempty stability condition, so that the unique 
unstable point in $\wt{\UU}^{ns}_{g,1}(g)$ is the origin, i.e., the point corresponding to the curve $C^{\cusp}(g)$.
We denote this point by $[C^{\cusp}(g)]$.
Then the functions $\a_{ij}[-m,q]$ identify the corresponding GIT quotient stack, 
$$\ov{\UU}^{ns}_{g,1}(g):=(\wt{\UU}^{ns}_{g,1}(g)\setminus [C^{\cusp}(g)])/\G_m,$$
with a closed substack in the weighted projective stack.

For two collection of weights as above, $\ba$ and $\ba'$, we denote by
$\wt{\UU}^{ns}_{g,n}(\ba,\ba')$ the interesection of the stacks $\wt{\UU}^{ns}_{g,n}(\ba)$ and $\wt{\UU}^{ns}_{g,n}(\ba')$.
In other words, we impose both conditions, $h^1(\sum a_ip_i)=0$ and $h^1(\sum a'_ip_i)=0$, on the marked points.

\subsection{The forgetful map}

The rational map \eqref{rat-map-eq} corresponds to a regular morphism
\begin{equation}\label{forg2-small-eq}
\forg_2:\wt{\UU}^{ns}_{g,2}((g-1,1),(g,0))\to \wt{\UU}^{ns}_{g,1}(g),
\end{equation}
which is given as the composition of the open embedding 
$\wt{\UU}^{ns}_{g,2}((g-1,1),(g,0))\hra \wt{\UU}^{ns}_{g,2}(g,0)$ followed by the forgetful map
$$\forg_2:\wt{\UU}^{ns}_{g,2}(g,0)\to \wt{\UU}^{ns}_{g,1}(g)$$
defined in \cite[Thm.\ A]{P-krich}. The latter map sends $(C,p_1,p_2,v_1,v_2)$, with $C$ irreducible, to $(C,p_1,v_1)$
(if $C$ is reducible then it gets replaced by a certain curve $\ov{C}$, such that $C\to \ov{C}$ is contraction of the component containing $p_2$).

Let $Z\sub \wt{\UU}^{ns}_{g,2}((g-1,1),(g,0))$ be the closed subscheme given as the preimage
of the origin under \eqref{forg2-small-eq}. Then there is a regular morphism
\begin{equation}\label{forg2-small2-eq}
\wt{\UU}^{ns}_{g,2}((g-1,1),(g,0))\setminus Z\to \ov{\UU}^{ns}_{g,1}(g)
\end{equation}
induced by \eqref{forg2-small-eq}.
Note that $Z$ consists of $(C,p_1,p_2,v_1,v_2)$ such that $(C,p_1)$ is the cuspidal curve $C^{\cusp}(g)$ (with the marked point
at infinity).

Let us denote by 
$$\wt{\UU}^{ns}_{g,2}((g-1,1),(g+1,0))\sub \wt{\UU}^{ns}_{g,2}(g-1,1)$$
the open subset given by the condition $h^1((g+1)p_1)=0$.
Let also 
$$\wt{\WW}\sub \wt{\UU}^{ns}_{g,2}((g-1,1),(g+1,0))$$ 
denote the closed locus given by the condition
$h^1(gp_1)\neq 0$, so that 
$$\wt{\UU}^{ns}_{g,2}((g-1,1),(g,0))=\wt{\UU}^{ns}_{g,2}((g-1,1),(g+1,0))\setminus \wt{\WW}.$$

Recall that we have sections $f_1[-m]\in H^0(C,\OO(mp_1+p_2))$, where $C$ is the universal curve
over $\wt{\UU}^{ns}_{g,2}(g-1,1)$, for $m\ge g$, with expansions at $p_1$ of
the form \eqref{f-i-m-eq} (with $i=1$) with $\a_1[-m,-g+1]=\a_1[-m,0]=0$.

\begin{lem}\label{ab-open-lem} 
Let us set $\a=\a_{12}[-g,-1]$, $\b=\a_{12}[-g-1,-1]$. Then the open subset 
$$\wt{\UU}^{ns}_{g,2}((g-1,1),(g,0))\sub \wt{\UU}^{ns}_{g,2}(g-1,1)$$ 
is given by the condition $\a\neq 0$.
Similarly, the open subset
$$\wt{\UU}^{ns}_{g,2}((g-1,1),(g+1,0))\sub \wt{\UU}^{ns}_{g,2}(g-1,1)$$ 
is the locus where either $\a\neq 0$ or $\b\neq 0$.
\end{lem}

\Pf .
Recall that the open subset $\wt{\UU}^{ns}_{g,2}((g-1,1),(g,0))$ 
is characterized by the condition $h^1(gp_1)=0$. Since $h^1(gp_1+p_2)=0$,
the long exact sequence of cohomology
associated with the exact sequence of sheaves
$$0\to \OO(gp_1)\to \OO(gp_1+p_2)\to \OO(p_2)/\OO\to 0$$
shows that $h^1(gp_1)\neq 0$ precisely for those curves for which $f_1[-g]$ is regular at $p_2$. But this is
equivalent to the vanishing of $\a$, since $\a$ is the coefficient of $t_2^{-1}$ in the expansion of $f_1[-g]$ at $p_2$.

The case of $\wt{\UU}^{ns}_{g,2}((g-1,1),(g+1,0))$ is similar: now we consider the exact sequence
$$0\to \OO((g+1)p_1)\to \OO((g+1)p_1+p_2)\to \OO(p_2)/\OO\to 0$$
which shows that $h^1((g+1)p_1)\neq 0$ when both $f_1[-g]$ and $f_1[-g-1]$ are regular at $p_2$,
i.e., both $\a$ and $\b$ vanish.
\ed

The following Proposition is a crucial step in proving Theorem A.

\begin{prop}\label{reg-map-prop} 
The subset $Z$ is closed in $\wt{\UU}^{ns}_{g,2}((g-1,1),(g+1,0))$, and we have
$Z\cap \wt{\WW}=\emptyset$. There exists a regular morphism
\begin{equation}\label{forg-2-main-eq}
\wt{\forg}_2:\wt{\UU}^{ns}_{g,2}((g-1,1),(g+1,0))\setminus Z\to \ov{\UU}^{ns}_{g,1}(g),
\end{equation}
extending the morphism \eqref{forg2-small2-eq} and
sending $\wt{\WW}$ to a point. Furthermore, this point has no nontrivial automorphisms.
\end{prop}

\Pf . Let $C'$ denote the universal curve over the open subset $\wt{\UU}^{ns}_{g,2}((g-1,1),(g,0))$.
To calculate explicitly the map \eqref{forg2-small-eq}, we need to find the sections
$f[-m]\in H^0(C',\OO(mp_1))$, for $m\ge g+1$, and a modified formal parameter $u$ at $p_1$, such that $f[-m]$
would have expansions of the form
\begin{equation}\label{f-m-u-eq}
f[-m]=u^{-m}+\a[-m,-g+1]u^{-g+1}+\a[-m,-g+2]u^{-g+2}+\ldots,
\end{equation}
where $\a[-m,q]$ are some rational expressions of the coordinates on $\wt{\UU}^{ns}_{g,2}(g-1,1)$ with only powers of
$\a$ in the denominator. 

As the first approximation let us set for $m\ge g+1$,
$$\wt{f}[-m]=f_1[-m]-\frac{\a_{12}[-m,-1]}{\a}f_1[-g].$$
The constant is chosen so that the poles at $p_2$ cancel out, so we have $\wt{f}[-m]\in H^0(\OO(mp_1))$,
while the expansion of $\wt{f}[-m]$ at $p_1$ has form
$$\wt{f}[-m]=t_1^{-m}-\frac{\a_{12}[-m,-1]}{\a}\cdot t_1^{-g}+\ldots,$$
where $t_1$ is the canonical parameter at $p_1$ on $C'$.

Now we need to change the canonical parameter to $u=t_1+c_1t_1^2+\ldots$, and to add to each $\wt{f}[-m]$ a
linear combination of $\wt{f}[-m']$ with $m'<m$, to get the expansions of the required form \eqref{f-m-u-eq}.
We want to know only the highest order polar parts of the funtions $\a[-m,q]$, i.e., those with the highest power of
$\a$ (prescribed below) in the denominator, so
we introduce the following filtration $F_n$ on the space of formal Laurent series in $t_1$ with coefficients in
$R=\OO(\wt{\UU}^{ns}_{g,2}((g-1,1),(g,0)))$. By definition, a Laurent series belongs to $F_n$ if it can be
written in the form $\sum_i a_i\a^{-i-n}t_1^i$, where each $a_i$ extends to a regular function on 
$\wt{\UU}^{ns}_{g,2}(g-1,1)$.

It will be enough for us to keep track only of $f[-m]\mod F_{m-1}$.
It is easy to see that the change of variables $t_1\mapsto t_1+c_1t_1^2+c_2t_1^3+\ldots$, where for each $i$, 
$\a^ic_i$ extends to a
regular function on $\wt{\UU}^{ns}_{g,2}(g-1,1)$, preserves the filtration $(F_n)$.
Since to go from $t_1$ to $u$ we will only use the changes of variables of this form, it suffices for us to know that
\begin{equation}\label{ti-f-g+1-mod-F-eq}
\wt{f}[-g-1]\equiv t_1^{-g-1}-\la t_1^{-g} \mod F_g,
\end{equation}
where $\la=\frac{\a_{12}[-g-1,-1]}{\a}=\frac{\b}{\a}$, while
\begin{equation}\label{ti-f-m-mod-F-eq}
\wt{f}[-m]\equiv t_1^{-m}\mod F_{m-1} \ \text{ for } m>g+1.
\end{equation}

We claim that there exist rational constants $(r_{m,j})$, $1\le j<m-g$, and $(r_i)$, $i\ge 1$, such that
\begin{equation}\label{f-la-formula-eq}
f[-m]\equiv \wt{f}[-m]+\sum_{1\le j<m-g}r_{m,j}\la^j \wt{f}[-m+j] \mod F_{m-1},
\end{equation}
for each $m\ge g+1$, and
\begin{equation}\label{t1-u-F-eq}
t_1\equiv u+r_1\la u^2+r_2\la^2 u^3+\ldots \mod F_{-2}.
\end{equation}
Namely, we prove by induction on $n\ge 1$ that \eqref{f-la-formula-eq} holds for all $m$ with $m\le g+n$, and
that the required relation between $t_1$ and $u$ holds modulo $t_1^nR[[t_1]]+F_{-2}$. 

Let us recall the recursive construction of $(f[-g-n])$ and of formal parameters $u_n$ such that 
$u_n\equiv u\mod t_1^{n+1}R[[t_1]]$, where $u$ is the canonical parameter (cf. \cite[Lem.\ 4.1.3]{FP}). 
For $n=1$ we have $f[-g-1]=\wt{f}[-g-1]$ and
$u_1=t_1$. Assume $f[-g-n']$ are already defined for $n'<n$ and $u_{n-1}\equiv u\mod t_1^nR[[t_1]]$ is known, 
so that 
$$f[-g-n']\equiv u_{n-1}^{-g-n'}\mod t_1^{-g+1}R[[t_1]] \text{ for } n'<n-1, \text{ while}$$
\begin{equation}\label{f-u-c-eq}
f[-g-n+1]\equiv u_{n-1}^{-g-n+1}+c\cdot u_{n-1}^{-g}\mod t_1^{-g+1}R[[t_1]].
\end{equation}
Then we set $u_n=u_{n-1}+\frac{c}{g+n-1}u_{n-1}^n$, the expansion of $f[-g-n+1]$ in $u_n$ will take form
$$f[-g-n+1]\equiv u_n^{-g-n+1} \mod t_1^{-g+1}R[[t_1]],$$
and the expansions of $f[-g-n']$ for all $n'<n-1$ in $u_n$ will still have the correct form. 
Now, if the expansion of $\wt{f}[-g-n]$ in $u_n$ has form
\begin{equation}\label{f-p-v-eq}
\wt{f}[-g-n]=u_n^{-g-n}+p_1u_n^{-g-n+1}+\ldots +p_{n-1}u_n^{-g-1}+\ldots,
\end{equation}
then we set
\begin{equation}\label{f-rec-formula-eq}
f[-g-n]=\wt{f}[-g-n]-p_1f[-g-n+1]-\ldots-p_{n-1}f[-g-1].
\end{equation}

The induction assumption implies that the function $c$ in \eqref{f-u-c-eq} has the leading polar term $r\la^{n-1}$ 
for some $r\in\Q$, so the change of variables from $u_{n-1}$ to $u_n$ is of the right form, as discussed above.
It follows that
$$t_1\equiv u_n+s_1\la u_n^2+\ldots+s_{n-1}\la^{n-1} u_n^n \mod t_1^{n+1}R[[t_1]]+F_{-2}$$
for some $s_i\in\Q$. Now from \eqref{ti-f-m-mod-F-eq} we get that
$$\wt{f}[-g-n]=(u_n+s_1\la u_n^2+\ldots+s_{n-1}\la^{n-1}u_n^n)^{-g-n} \mod t_1^{-g}R[[t_1]]+F_{g+n-1}.$$ 
This implies that for $i=1,\ldots,n=1$, the leading polar term of the coefficient $p_i$ in the expansion
\eqref{f-p-v-eq} is of the form $a_i\la^i$, for $a_i\in\Q$. Now \eqref{f-rec-formula-eq} shows that
\eqref{f-la-formula-eq} holds for $m=g+n$.
This finishes the proof of our claim.

Now combining \eqref{ti-f-g+1-mod-F-eq}--\eqref{t1-u-F-eq}, we get that for each $m\ge g+1$ the expansion
of $f[-m]$ in the canonical parameter $u$ has form 
$$f[-m]\equiv u^{-m}+\sum_{j\ge 1} s_{m,j}\la^{m-g+j}u^{-g+j}\mod F_{m-1},$$
for some rational constants $(s_{m,j})$. In other words, the functions $\a[-m,-g+j]\in R$, defining the map
\eqref{forg2-small-eq}, have form
$$\a[-m,-g+j]=s_{m,j}\la^{m-g+j}+\ldots$$
where the omitted terms have smaller powers of $\a$ in the denominator.

Finally, we need to know that not all $(s_{m,j})$ are zero, so let us compute $s_{-g-1,-g+1}$ and
$s_{g-1,-g+2}$ following the above procedure (we will need to look at two coordinates to prove that the point, which is the image of $\WW$, has no nontrivial automorphisms).
Due to \eqref{ti-f-g+1-mod-F-eq}, the first change of variables is 
$$t_1=u_2-\frac{\la}{g+1}u_2^2\mod F_{-2}.$$
Then we get expansions
\begin{align*}
&f[-g-1]=\wt{f}[-g-1]\equiv u_2^{-g-1}+\frac{2-g}{2(g+1)}\la^2 u_2^{-g+1}+\\
&\frac{-g^2+g+3}{3(g+1)^2}\la^3 u_2^{-g+2} \ \mod u_2^{-g+3}R[[u_2]]+F_g,
\end{align*}
$$\wt{f}[-g-2]\equiv u_2^{-g-2}+\frac{g+2}{g+1}\la u_2^{-g-1}+\frac{(g+2)(g+3)}{2(g+1)^2}\la^2 u_2^{-g} \ \mod u_2^{-g+1}R[[u_2]]
+F_{g+1},$$
\begin{align*}
&\wt{f}[-g-3]\equiv u_2^{-g-3}+\frac{g+3}{g+1}\la u_2^{-g-2}+\frac{(g+3)(g+4)}{2(g+1)^2}\la^2 u_2^{-g-1}+\\
&\frac{(g+3)(g+4)(g+5)}{6(g+1)^3}\la^3 u_2^{-g} \ \mod u_2^{-g+1}R[[u_2]]
+F_{g+2}.
\end{align*}
Hence, the coefficient of $u_2^{-g}$ in $f[-g-2] \mod F_{g+1}$ (which is the same as in $\wt{f}[-g-2] \mod F_{g+1}$) is 
$\frac{(g+2)(g+3)}{2(g+1)^2}\la^2$.
Thus, the second change of variables (defined so that the coefficient of $u_3^{-g}$ in $f[-g-2]$ is zero) is 
$$u_2=u_3+\frac{(g+2)(g+3)}{2(g+2)(g+1)^2}\la^2 u_3^3 \mod F_{-2},$$
and we get the expansion
$$f[-g-1]=u_3^{-g-1}-\frac{2g+1}{2(g+1)}\la^2 u_3^{-g+1}+\frac{-g^2+g+3}{3(g+1)^2}\la^3 u_3^{-g+2} \ \ \mod u_3^{-g+3}R[[u_3]]+F_g,$$
which shows that 
$$s_{g+1,1}=-\frac{2g+1}{2(g+1)}.$$ 
Also, we see that the coefficient of $u_3^{-g}$ in the expansion
of $\wt{f}[-g-3] \mod F_{g+2}$ is equal to $-\frac{(g+3)(g^2+3g-1)}{3(g+1)^3}\la^3$. This dictates that the next change of variables
is 
$$u_3=u_4-\frac{(g^2+3g-1)}{3(g+1)^3}\la^3 u_4^4 \mod F_{-2}.$$
Finally, we get that the coefficient of $u_4^{-g+2}$ in the expansion of $f[-g-1] \mod F_g$ is equal
to 
$$\frac{-g^2+g+3}{3(g+1)^2}\la^3+\frac{(g^2+3g-1)}{3(g+1)^2}\la^3=\frac{4g+2}{3(g+1)^2}\la^3,$$
and hence, 
$$s_{g+1,2}=\frac{4g+2}{3(g+1)^2}.$$

Now let us consider the modified map
$$\a\cdot\wt{\forg}_2: \wt{\UU}^{ns}_{g,2}((g-1,1),(g,0))\to \wt{\UU}^{ns}_{g,1}(g): x\mapsto \a(x)\cdot \wt{\forg}_2(x)$$
Since the weight of $\a[-m,-g+j]$ is $m-g+j$, the modifed map sends $x$ to the point in 
$\wt{\UU}^{ns}_{g,1}(g)$ with coordinates
\begin{equation}\label{a-rescaled-coord-eq}
\a(x)^{m-g+j}\a[-m,-g+j](x)=s_{m,j}\b(x)^{m-g+j}+\a(x)\cdot f_{m,j}(x),
\end{equation}
where $f_{m,j}$ are regular functions on $\wt{\UU}^{ns}_{g,2}(g-1,1)$.
In particular, $\a\cdot\wt{\forg}_2$ can be viewed as a regular map from $\wt{\UU}^{ns}_{g,2}(g-1,1)$.

Recall that by Lemma \ref{ab-open-lem}, the open subset $\wt{\UU}^{ns}_{g,2}((g-1,1),(g+1,0))$ is the locus where
either $\a\neq 0$ or $\b\neq 0$, and the locus $\wt{\WW}$ is given by $\a=0$.
Thus, \eqref{a-rescaled-coord-eq} gives for $x\in\wt{\WW}$:
$$\a\cdot\wt{\forg}_2(x)=(s_{m,j}\b(x)^{m-g+j})=\b(x)\cdot (s_{m,j}).$$ 
Furthermore, as we have seen above, the constants $s_{-g-1,1}$ and
$s_{-g-1,2}$ are nonzero, so the corresponding coordinates in the above expression are also nonzero.
Note also that the corresponding point of $\ov{\UU}^{ns}_{g,1}(g)$ is equal to $(s_{m,j})$, 
so it does not depend on $x$. 

Denoting by $U_{\b\neq 0}\sub \wt{\UU}^{ns}_{g,2}((g-1,1),(g+1,0))$ the open subset where $\b\neq 0$,
we get
$$(\a\cdot\wt{\forg}_2)^{-1}(0)\cap U_{\b\neq 0}=Z\cap U_{\b\neq 0},$$
and so $Z\cap U_{\b\neq 0}$ is closed in $U_{\b\neq 0}$. Since, $Z$ is closed in the open subset $\a\neq 0$, we derive
that $Z$ is closed in $\wt{\UU}^{ns}_{g,2}((g-1,1),(g+1,0))$.

We have a covering of $\wt{\UU}^{ns}_{g,2}((g-1,1),(g+1,0))\setminus Z$ by two open subsets:
$\wt{\UU}^{ns}_{g,2}((g-1,1),(g,0))\setminus Z$ and $U_{\b\neq 0}$. The required regular morphism \eqref{forg-2-main-eq}
to $\ov{\UU}^{ns}_{g,1}(g)$ is induced by $\wt{\forg}_2$ on $\wt{\UU}^{ns}_{g,2}((g-1,1),(g,0))\setminus Z$ and 
by $\a\cdot\wt{\forg}_2$ on $U_{\b\neq 0}$. As we have seen above, this morphism sends
$\wt{\WW}\sub U_{\b\neq 0}$ to the point $(s_{m,j})$ of the weighted projective stack
with two nonzero homogeneous coordinates, of weights $2$ and $3$.
Hence, this point does not have nontrivial automorphisms.
\ed

\subsection{Proof of Theorem A}


It is well known that $\WW\sub \MM_{g,1}$ is an irreducible divisor (see \cite{Arb1,Arb2}).
Now let $U\sub\MM_{g,1}$ be the open substack of $(C,p)$ satisfying 
$$h^1((g+1)p)=0, \ \ h^0((g-1)p)=1.$$
Note also that we have an inclusion
$$\MM_{g,1}\setminus \WW\sub U$$
since the condition $h^0(gp)=1$ implies that $h^1((g+1)p)=0$ and $h^0((g-1)p)=1$.
Furthermore, the complement to $U$ is a proper closed subset in $\WW$, so it has codimension $\ge 2$ in $\MM_{g,1}$. 
In particular, $U\cap \WW$ is dense in $\WW$. 

Note that we have a natural open inclusion 
\begin{equation}\label{reg-map-eq}
\MM_{g,1}\setminus \WW\hra \ov{\UU}^{ns}_{g,1}(g).
\end{equation}
Indeed, the only unstable point in $\UU^{ns}_{g,1}(g)$ corresponds to the singular curve $C^{\cusp}(g)$.
We are going to show that the above morphism extends to a regular morphism
$$U\to \ov{\UU}^{ns}_{g,1}(g),$$
such that $U\cap \WW$ is mapped to a point.

Recall that by Proposition \ref{reg-map-prop}, we have a regular morphism 
$$\wt{\forg}_2:\wt{\UU}^{ns}_{g,2}((g-1,1),(g+1,0))\setminus Z\to \ov{\UU}^{ns}_{g,1}(g),$$
sending $\wt{\WW}$ to a point. Let 
$V\sub \wt{\UU}^{ns}_{g,2}((g-1,1),(g+1,0))$
be the open subset corresponding to smooth curves. Then $V\cap Z=\emptyset$ because for points of $Z$
the underlying curve is singular. Thus, the above morphism induces a regular morphism
\begin{equation}\label{forg2-V-eq}
\wt{\phi}:V\to \ov{\UU}^{ns}_{g,1}(g),
\end{equation}
mapping $\wt{\WW}\cap V$ to a point.

Now we claim that the natural projection $V\to \MM_{g,1}$ induces a smooth surjective morphism $V\to U$.
Indeed, if $h^0((g-1)p_1+p_2)=1$ then $h^0((g-1)p_1)=1$, so this projection factors through $U$.
Conversely, if for $(C,p_1)\in\MM_{g,1}$ one has $h^0((g-1)p_1)=1$ then for generic $p_2$ we will have $h^0((g-1)p_1+p_2)=1$,
hence the map $V\to U$ is surjective. It is smooth since $V$ is a $\G_m^2$-torsor over an open
substack of a universal curve over $U$.

It remains to prove that the morphism \eqref{forg2-V-eq} factors through a morphism $\phi:U\to \ov{\UU}^{ns}_{g,1}(g)$
(it will then map $\WW\cap U$ to a point, since \eqref{forg2-V-eq} sends $\wt{\WW}\cap V$ to a point).
Indeed, this is true if we restrict to the open subset $\MM_{g,1}\setminus\WW$, by the construction.
Now let us set $T:=V\times_U V$ and consider two morphisms 
$$f_1=\wt{\phi}\circ \pi_1,f_2=\wt{\phi}\circ \pi_2: T\to \ov{\UU}^{ns}_{g,1}(g),$$
where $\pi_1$ and $\pi_2$ are two projections to $V$. We know that these two maps agree on the open
subset $\pi^{-1}(\MM_{g,1}\setminus\WW)$, where $\pi$ is the projection $T=V\times_U V\to U$.

Note that the scheme $T$ parametrizes data $(C,p_1,p_2,p'_2,v_1,v_2,v'_2)$ such that 
$h^0((g-1)p_1+p_2)=h^0((g-1)p_1+p'_2)=1$ and $h^1((g+1)p_1)=0$ (and $C$ smooth, $p_1\neq p_2$, $p_1\neq p'_2$).
Thus, it is an open subset in a $\G_m^3$-torsor over the universal curve over $\MM_{g,2}$ 
(via the projection to $(C,p_1,p_2,p'_2)$), in particular, $T$ is smooth and irreducible. 

Let us consider the cartesian diagram
\begin{diagram}
T'&\rTo{}&\ov{\UU}^{ns}_{g,1}(g)\\
\dTo{\rho}&&\dTo{\De}\\
T&\rTo{(f_1,f_2)}&\ov{\UU}^{ns}_{g,1}(g)\times \ov{\UU}^{ns}_{g,1}(g)
\end{diagram}
Since the stack $\ov{\UU}^{ns}_{g,1}(g)$ is separated, the vertical arrows are finite morphisms.
Finally, we observe that a generic pointed curve $(C,p)$ in $\MM_{g,1}$ does not have nontrivial automorphisms (note that in
the case $g=2$ this is true since we can take $p$ not to be a Weierstrass point). Hence,
the preimages of points with trivial automorphisms in $\ov{\UU}^{ns}_{g,1}(g)$ under $f_1$ and $f_2$ are nonempty
open subsets in $T$. Since $f_1$ and $f_2$ agree on a nonempty open subset, we deduce
that there exists a nonempty open subset $W\sub T$ such that $\rho^{-1}(W)\to W$ is an isomorphism.
Let $T''\sub T'$ be an irreducible component of $T'$, containing $\rho^{-1}(W)$, with reduced scheme structure.
Then $\rho|_{T''}:T''\to T$ is a finite birational morphism. Since $T$ is smooth, we deduce that $\rho|_{T''}$ is an isomorphism.
Hence, $\rho$ admits a section, and so we have $f_1=f_2$, which means that the map 
\eqref{forg2-V-eq} descends to a morphism from $U$.
\ed

\section{Curves of genus $2$}\label{gen-2-sec}

\subsection{Explicit identification of $\wt{\UU}^{ns}_{2,1}(2)$}\label{g-2-n-1-a-2-sec}

\begin{prop}\label{1pt-moduli-prop} 
Let us work over $\Z[1/6]$.
One has an isomorphism of the moduli scheme
$\wt{\UU}^{ns}_{2,1}(2)$ with the affine space $\A^5$ with coordinates $q_1, q_{2,0}, q_{2,1}, q_{3,0}, q_{3,1}$,
so that the affine universal curve $C\setminus \{p\}$ is given by the following equations in the independent variables $f,h,k$:
\begin{equation}\label{g-2-n-1-eq}
\begin{array}{l}
h^2=fk+q_1h+2q_1^2+f(q_{2,0}+q_{2,1}f),\\
hk=f(q_{3,0}+q_{3,1}f+f^2)-q_1k+(q_{2,0}+q_{2,1}f)h+q_1(q_{2,0}+q_{2,1}f),\\
k^2=(q_{3,0}+q_{3,1}f+f^2)h+(q_{2,0}+q_{2,1}f)^2-2q_1(q_{3,0}+q_{3,1}f+f^2).
\end{array}
\end{equation}
The weights of the $\G_m$-action are:
$$\deg(q_{2,1})=2, \ \deg(q_{3,1})=3, \ \deg(q_1)=4, \ \deg(q_{2,0})=5, \ \deg(q_{3,0})=6.$$
Hence, we get the identification of $\ov{\UU}^{ns}_{2,1}(2)$ with the weighted projective stack
$\P(2,3,4,5,6)$.
\end{prop}



\Pf . This is proved using the same method as in \cite[Thm.\ A]{P-ainf} and \cite[Thm.\ A]{P-krich}.
Let $(C,p,v)$ be a point in $\wt{\UU}^{ns}_{2,1}(2)$. Since $h^1(2p)=0$, we have $h^0(np)=n-1$ for $n\ge 2$.
Let $t$ be a formal parameter at $t$ compatible with the given tangent vector.
We can find the elements $f\in H^0(C,\OO(3p))$, $h\in H^0(C,\OO(4p))$ and $k\in H^0(C,\OO(5p))$ with the Laurent expansions
$$f=\frac{1}{t^3}+\ldots, \ \ h=\frac{1}{t^4}+\ldots, \ \ k=\frac{1}{t^5}+\ldots,$$
where the omitted terms have poles of smaller order. Then the elements 
\begin{equation}\label{algebra-lin-basis-eq}
f^n, f^nh, f^nk, \ \text{ for } n\ge 0,
\end{equation}
form a linear basis on 
$H^0(C\setminus\{p\},\OO)$, so we can express $h^2$, $hk$ and $k^2$ as their linear combinations.
Taking into account the above Laurent expansion, we get relations of the form
\begin{equation}\label{g-2-n-1-prelim-eq}
\begin{array}{l}
h^2=p_1(f)k+q_1(f)h+c_1(f),\\
hk=p_2(f)k+q_2(f)h+c_2(f),\\
k^2=p_3(f)k+q_3(f)h+c_3(f),
\end{array}
\end{equation}
where $p_i$, $q_i$, $c_i$ are polynomials in $f$ with the following restrictions: 
\begin{align*}
&\deg p_1=1, \deg p_2\le 1, \deg p_3\le 1, \deg q_1\le 1, \deg q_2\le 1, \deg q_3=2, \\
&\deg c_1\le 2, \deg c_2=3,
\deg c_3\le 3,
\end{align*}
and the polynomials $p_1$, $q_3$ and $c_2$ are monic.
Note that $f$ is defined up to adding a constant, while  $h$ and $k$ are defined up to the transformation
$$(h,k)\mapsto (\wt{h}=h+A(f), \wt{k}=k+Bh+C(f)),$$
where $A$ and $C$ are linear polynomials in $f$ and $B$ is a constant.
It is easy to check that we can fix the ambiguity in the choice of $h$ and $k$ by requiring that $p_3=0$
and $p_2=-q_1$ is a constant, i.e., does not have a linear term in $f$. 
More precisely, we should set
\begin{equation}
A=-\frac{q_1+p_2}{3}, \ \ B=\frac{1}{3}(q'_1-2p'_2), \ \ C=-\frac{p_3}{2}-\frac{B^2p_1}{2}-Bp_2
\end{equation}
(here $q'_1$ and $p'_2$ are derivatives of the linear polynomials $q_1$ and $p_2$).
Note that here we use our assumption that $6$ is invertible. 
Finally, we can fix the ambiguity in the choice of $f$ by requiring that $p_1(f)=f$.

Now the fact that the elements \eqref{algebra-lin-basis-eq} form a basis of $H^0(C\setminus\{p\},\OO)$ is equivalent
to the condition that the relations \eqref{g-2-n-1-prelim-eq} form a Gr\"obner basis in the ideal they generate
(with respect to the degree reverse lexicographical order such that $f<h<k$, $\deg(f)=3$, $\deg(h)=4$, $\deg(k)=5$).
Applying the Buchberger's Criterion (see \cite[Thm.\ 15.8]{Eis-CA}) we compute that this condition 
is equivalent to the following expressions of $c_1,c_2,c_3$ in terms of the other variables 
(where in the second expression in each line we take into account the normalization $p_3=0$, $p_2=-q_1$): 
\begin{align*}
&c_1=p_2^2+ p_1q_2-q_1p_2=2q_1^2+p_1q_2,\\
&c_2=p_1q_3-p_2q_2=p_1q_3+q_1q_2,\\
&c_3=q_2^2+p_2q_3-q_1q_3=q_2^2-2q_1q_3.
\end{align*}
Thus, if we set 
$$q_2=q_{2,0}+q_{2,1}f, \ \ q_3=q_{3,0}+q_{3,1}f+f^2,$$
then we see that the constants $(q_1,q_{2,0}, q_{2,1}, q_{3,0}, q_{3,1})$ determine the curve $(C,p)$.
The above process can be run in families and can be reversed (see the proofs of \cite[Thm.\ A]{P-ainf} and \cite[Thm.\ A]{P-krich}),
so this gives the required identification of our moduli space with $\A^5$.
\ed



\subsection{Special cuspidal curve $C_0$}\label{cusp-sec}

Let $C_0$ denote the curve obtained from $\P^1$ by pinching the point $0$ into a genus $2$ cuspidal singular point,
so that a regular function $f$ near $0$ descends to $C_0$ if and only if the expansion of $f$ in the standard parameter $t$
has form 
\begin{equation}\label{C0-exp-def}
f\equiv c_0+c_2\cdot t^2\mod (t^4).
\end{equation}
Note that this condition depends on coordinates, i.e., the point $\infty\in C_0$ plays a special role.
For example, the standard $\G_m$-action on $\P^1$, preserving $0$ and $\infty$, descends to a $\G_m$-action on $C_0$.
Also, note that $C_0\setminus\{\infty\}=\Spec(\C[t^2,t^5])$.

The next Lemma shows that 
if we equip $C_0$ with a smooth marked point $p\neq\infty$ then we get a point of $\ov{\UU}^{ns}_{2,1}(2)$.

\begin{lem}\label{C0-lem} 
Let $p\in C_0\setminus \{0,\infty\}$. Then $h^0(C_0,\OO(2p))=1$. On the other hand, for $p=\infty$ we
have $h^0(C_0,\OO(2p))=2$.
\end{lem}

\Pf . In the case $p\neq 0,\infty$ we can assume that $t(p)=1$. 
Then $\OO_{\P^1}(2p)$ is spanned by $1$, $\frac{1}{1-t}$ and $\frac{1}{(1-t)^2}$.
Looking at the expansions at $t=0$ we see that the only sections of $\OO_{\P^1}(2p)$ satisfying \eqref{C0-exp-def}
are constants.

In the case $p=\infty$ the functions $(1, t^2)$ give a basis of $H^0(C_0,\OO(2p))$. 
\ed

\begin{defi}\label{C0-defi}
We denote by $[C_0]$ the point of $\ov{\UU}^{ns}_{2,1}(2)$ corresponding to $(C_0,p)$, where $p\neq 0,\infty$.
\end{defi}

\subsection{Classification of singular irreducible curves of genus $2$}\label{classif-sec}

Let $C$ be an irreducible curve of genus $2$, and let $\rho:\wt{C}\to C$ be the normalization. If $C$ is singular then
the genus of $\wt{C}$ is either $1$ or $0$. 

If the genus of $\wt{C}$ is $1$ then $\coker(\OO_C\to \rho_*\OO_{\wt{C}})$
has length $1$, so it is supported at one singular point $q\in C$. If $\rho^{-1}(q)$ contains two distinct points $q_1,q_2\in C$
then $\rho$ factors through a morphism $C'\to C$, where $C'$ is the nodal curve
obtained by gluing $q_1$ and $q_2$ on $\wt{C}$. Since $C'$ has genus $2$ we should have $C\simeq C'$.
If $\rho^{-1}(q)$ is one point on $C$ then it is easy to see that $C$ has a simple cusp at $q$.

In the remaining case when $\wt{C}=\P^1$ we have more possibilities. The length of the sheaf 
$\FF:=\coker(\OO_C\to \rho_*\OO_{\wt{C}})$
is now $2$, so the support of $\FF$ can consist of $\le 2$ points. 

\medskip

\noindent 
{\bf Case I: support of $\FF$ consists of two distinct points $q_1$, $q_2$}.
We have the following subcases.

{\bf Case Ia: $|\rho^{-1}(q_1)|>1$ and $|\rho^{-1}(q_2)|>1$}. 
In this case the map $\rho$ factors through the nodal curve $C'$ obtained
by gluing two pairs of distinct points in $\P^1$. Since the genus of $C'$ is $2$, we should have $C\simeq C'$.

{\bf Case Ib: $|\rho^{-1}(q_1)|=1$ and $|\rho^{-1}(q_2)|>1$}. In this case $\rho$ factors through the 
curve $C'$ obtained by gluing two
pairs of distinct point in $\P^1$ and pinching one extra point to a simple cusp. Again, we have that the genus of $C'$ is $2$,
so $C\simeq C'$.

{\bf Case Ic: $|\rho^{-1}(q_1)|=|\rho^{-1}(q_2)|=1$}. In this case $C$ is obtained by pinching two points of $\P^1$ into simple
cusps.

\medskip

\noindent
{\bf Case II: $\FF$ is supported at one point $q$}.

{\bf Case IIa: $|\rho^{-1}(q)|>2$}. In this case $\rho$ factors through the curve $C'$ obtained by gluing transversally $3$
points on $\P^1$ into a single point (with the coordinate cross singularity). Since the genus of $C'$ is $2$, we get
$C\simeq C'$.

{\bf Case IIb: $|\rho^{-1}(q)|=2$}. Let $\rho^{-1}(q)=\{q_1,q_2\}$. Let $t$ be a generator of the maximal ideal 
$\mg_q\sub \OO_{C,q}$. Assume first that $t\in \mg_{q_1}^2$. 
Then $\rho$ factors through the curve $C'$ obtained from $\P^1$ by first
pinching $q_1$ into a simple cusp and then gluing it transversally with the point $q_2$. Since $C'$ has genus $2$, we have
$C\simeq C'$. On the other hand,
if $t$ maps to a generator of $\mg_{q_i}$ for $i=1,2$, then $\rho$ factors through the curve $C'$ obtained
from $\P^1$ by gluing $q_1$ and $q_2$ into a tacnode singularity. Since such $C'$ has genus $2$, we have $C\simeq C'$.

{\bf Case IIc: $|\rho^{-1}(q)|=1$}. In this case we can identify $C$ with $\P^1$ as a topological space, so that
$\OO_C$ is a subsheaf of $\OO_{\P^1}$, which differs from it only at one point $q$, so that $\mg_{C,q}\sub \mg^2_{\P^1,q}$
is an embedding of codimension $1$. We claim that there are two curves of this type, up to an isomorphism.
If $\mg_{C,q}\sub \mg^3_{\P^1,q}$ then $\mg_{C,q}=\mg^3_{\P^1,q}$ and
$C=C^{\cusp}(2)$ (see Sec.\ \ref{background-sec}). Now assume that $\mg_{C,q}\not\sub \mg^3_{\P^1,q}$. 
Let $t$ be a formal parameter near $q$ on $\P^1$. Then $\hat{\mg}_{C,q}$ is a (non-unital) subalgebra in $t^2\C[[t]]$ of
codimension $1$, and there exists an element $f\in\hat{\mg}_{C,q}$ such that 
$f\equiv t^2\mod t^3\C[[t]]$. Changing the formal parameter we can assume that $f=t^2$.
There could not be an element $h\in\hat{\mg}_{C,q}$ such that 
$h\equiv t^3\mod t^4\C[[t]]$, since then we would have $\hat{\mg}_{C,q}=t^2\C[[t]]$. Therefore, 
$$\hat{\mg}_{C,q}=\C\cdot t^2+t^4\C[[t]].$$
Note that the subspace in the right-hand side depends only on $t\mod t^3\C[[t]]$.
Now we observe that any formal parameter at $q$, modulo $\mg^3_{\P^1,q}$, can be obtained from a unique
regular function on $\P^1\setminus\{p\}$, for some $p\neq q$. Using automorphisms of $\P^1$ we can make $q=0$,
$p=\infty$, so that $C$ is the curve $C_0$ defined before.

\subsection{Comparison of stabilities for irreducible curves of genus $2$}

\begin{prop}\label{g2-Zstable-prop} 
Let $C$ be an irreducible curve of genus $2$, and let $p$ be a smooth point. Then $(C,p)$
$\ZZ$-stable if and only if $C$ is not of type IIc.
\end{prop}

\Pf . It is easy to see that a curve $C$ of type IIc is not $\ZZ$-stable. Indeed, if there is a contracting map $C'\to C$
then $C'$ would have a rational component with only two distinguished points, so it could not be stable.
Assume now that $C$ is not of type IIc.
If $(C,p)$ is nodal then it is stable (since $C$ is irreducible), hence it is $\ZZ$-stable.

Next, if $C$ is obtained by pinching a point on an irreducible nodal curve $E$ of genus $1$ into a cusp, then there is a contraction
$f:E\cup E'\to C$, where $E\cup E'$ is the stable curve with $E$ and $E'$ glued nodally at one point.
Here the marked point is placed on $E$ and $f(E')$ is the cusp on $C$. This shows that $(C,p)$ is $\ZZ$-stable.
Similarly, if $C$ is a rational curve with two cusps then there is a contraction to $C$ from $\P^1$ with two elliptic tails
(that get contracted into cusps). 

There remains two cases for $C$: IIa and IIb. In the case IIa we have a contraction to $C$ from the union of two
$\P^1$'s, joined nodally at $2$ points. In the case IIb there is a contraction to $C$ from the curve with an elliptic bridge.
In other words, we consider the union $\P^1\cup E$, where $E$ is an elliptic curve, $\P^1$ and $E$ are 
joined nodally at $2$ points, so that there are
no marked points on $E$. It is known (see \cite[Ex.\ 2.5]{Smyth-modular}) 
that there exists a contraction $\P^1\cup E\to C$, mapping $E$ to the singular point, 
for both types of curves occurring in the case IIb.
\ed

\begin{cor}\label{smooth-cor}
The stack $\ov{\MM}_{2,1}(\ZZ)$ is smooth and irreducible.
\end{cor}

\Pf . The possible singular points that can appear in $\ZZ$-stable curves of genus $2$, other that nodes, are:
a simple cusp, a tacnode, and a coordinate cross in $3$-space. 
All of these have smooth versal deformation spaces and are smoothable,
hence the assertion (see \cite[Lem.\ 2.1]{Smyth-II}).
\ed

Using the classification from Sec.\ \ref{classif-sec} we easily get the following 
codimension estimate.

\begin{lem}\label{g2-non-nodal-lem}
Away from a closed subset of codimension $\ge 2$, for 
every point $(C,p)$ in $\ov{\MM}_{2,1}(\ZZ)$ (resp., $\ov{\UU}^{ns}_{2,1}(2)$),
$C$ is either smooth, or a nodal curve with the normalization of genus $1$.
\end{lem}

\Pf . Both $\ov{\MM}_{2,1}(\ZZ)$ and $\ov{\UU}^{ns}_{2,1}(2)$ are irreducible of dimension $4$.
Now we just go through the strata described in Sec.\ \ref{classif-sec}
and see that they all have dimension
$\le 2$, except when $C$ is either smooth or nodal with the normalization of genus $1$.
\ed

We need one more simple observation.

\begin{lem}\label{g2-ab-open-set-lem} 
Let $C$ be an irreducible curve of genus $2$, and let $p\in C$ be a smooth point. Then $h^0(p)=1$ and
$h^1(3p)=0$.
\end{lem}

\Pf . First, if $h^0(p)=2$ then we would get a degree $1$ regular map $C\to\P^1$. Composing
it with the normalization map $\wt{C}\to C$, we get that the normalization map is the inverse map $\P^1\to C$,
which is impossible. Hence, $h^0(p)=1$.

If $h^1(2p)=0$ then we also have $h^1(3p)=0$, so it is enough to consider the case $h^1(2p)\neq 0$, i.e., $h^0(2p)=2$.
Suppose that $h^0(3p)=3$. Then we can choose $f\in H^0(C,\OO(2p))$ and $h\in H^0(C,\OO(3p))$ with the Laurent
expansions $f=\frac{1}{t^2}+\ldots$, $h=\frac{1}{t^3}+\ldots$ at $p$ (for some formal parameter $t$ at $p$).
Furthermore, there is a canonical choice of $f$ and $h$, such that the relation
$$h^2=f^3+af+b$$
holds for some constants $a$ and $b$. 
Then the algebra $\OO(C\setminus \{p\})$ has the linear basis $(f^n)$, $(hf^n)$, and is isomorphic
to the algebra $A=\C[h,f]/(h^2-f^3-af-b)$. Since $C$ is irreducible, it is isomorphic to $\Proj$ of the Rees algebra of $A$,
which is a plane cubic, so we get that the arithmetic genus of $C$ is equal to $1$, which is a contradiction.
This shows that $h^0(3p)=2$, i.e., $h^1(3p)=0$.
\ed

\begin{thm}\label{g2-thm}
Let $\ov{\WW}\sub \ov{\MM}_{2,1}(\ZZ)$ be the closure of the Weierstrass locus $\WW\sub \MM_{2,1}$.
Then $\ov{\WW}$ coincides with the locus where $h^1(2p)\neq 0$.
There is a regular morphism
$$\phi_2:\ov{\MM}_{2,1}(\ZZ)\to \ov{\UU}^{ns}_{2,1}(2),$$
such that $\phi_2(\ov{\WW})=[C_0]$ and
$\phi_2$ induces an isomorphism
$$\ov{\MM}_{2,1}(\ZZ)\setminus \ov{\WW}\rTo{\sim} \ov{\UU}^{ns}_{2,1}(2)\setminus [C_0].$$
\end{thm}

\Pf . First, we observe that every irreducible component of the locus $h^1(2p)\neq 0$ has codimension $1$ in 
$\ov{\MM}_{2,1}(\ZZ)$ (recall that the latter stack is smooth and irreducible by Corollary \ref{smooth-cor}).
By Lemma \ref{g2-non-nodal-lem}, to see that this locus coincides with $\ov{\WW}$, it is enough to see
that the locus of $(C,p)$, such that $C$ is nodal with normalization $E$ of genus $1$ and $h^1(2p)\neq 0$ has
dimension $2$ (and hence has codimension $2$ in $\ov{\MM}_{2,1}(\ZZ)$). But if $C$ is obtained from $E$ by identifying
points $q_1\neq q_2$ then the condition that $h^0(2p)=2$ implies the existence of a rational function on $E$ with pole
of order $2$ at $p$ and vanishing at both $q_1$ and $q_2$. In other words, we should have a linear equivalence
$2p\sim q_1+q_2$. Thus, we have a finite number of choices fo each $(E,q_1,q_2)$, so the dimension is $2$.

Next, let us denote by
$$V^{\ZZ}\sub \wt{\UU}^{ns}_{2,2}(1,1)$$
the open substack consisting of $(C,p_1,p_2)$ such that $(C,p_1)$ is $\ZZ$-stable (in particular, $C$ is irreducible).
Lemma \ref{g2-ab-open-set-lem} shows that every $(C,p)$ in $\ov{\MM}_{2,1}(\ZZ)$ satisfies
$h^0(p)=1$ and $h^1(3p)=0$. This implies that
$$V^{\ZZ}\sub \wt{\UU}^{ns}_{2,2}((1,1),(3,0))$$
and the projection $V^\ZZ\to \ov{\MM}_{2,1}(\ZZ)$ is surjective.
Furthermore, since the curve $[C^{\cusp}(2)]$ is not $\ZZ$-stable, we have the inclusion
$$V^{\ZZ}\sub \wt{\UU}^{ns}_{2,2}((1,1),(3,0))\setminus Z.$$
Thus, the restriction of the map \eqref{forg-2-main-eq} gives us a regular morphism
\begin{equation}\label{Zstable-forg-map}
V^\ZZ\to \ov{\UU}^{ns}_{2,1}(2),
\end{equation}
contracting $\wt{\WW}$ to a point. 

Now, similarly to the proof of Theorem A we check that the morphism \eqref{Zstable-forg-map} factors through 
$\ov{\MM}_{2,1}(\ZZ)$. Note that to apply the same argument as in Theorem A we use the following facts: 
(i) $\ov{\MM}_{2,1}(\ZZ)$ is smooth (see Corollary \ref{smooth-cor}); (ii) the projection
$V^\ZZ\to \ov{\MM}_{2,1}(\ZZ)$ is smooth (since $p_2$ varies in a smooth part of a curve); and
(iii) $V^\ZZ\times_{\ov{\MM}_{2,1}(\ZZ)} V^\ZZ$ is irreducible, as a $\G_m^3$-torsor over the moduli
stack of $(C,p_1,p_2,p'_2)$ with $C$ smoothable.

This gives us the required morphism $\phi_2$ contracting to $\ov{\WW}$ to some point in $\ov{\UU}^{ns}_{2,1}(2)$.
On the other hand, by Proposition \ref{g2-Zstable-prop}, the only point in $\ov{\UU}^{ns}_{2,1}(2)$, which
is not $\ZZ$-stable is $[C_0]$ (recall that by this we mean the pointed curve $(C_0,p)$, where $p\neq 0,\infty$, see
Lemma \ref{C0-lem}). 
Thus, the rational map $\phi_2^{-1}$ is regular on $\ov{\UU}^{ns}_{2,1}(2)\setminus [C_0]$ (and sends $(C,p)$ to $(C,p)$).
Also, the restriction of $\phi_2$ to $\ov{\MM}_{2,1}(\ZZ)\setminus\ov{\WW}$, i.e., to the locus where $h^0(2p)=1$,
is an open embedding sending $(C,p)$ to $(C,p)$.
This implies that $\phi_2(\ov{\WW})=[C_0]$, and $\phi_2$ induces an isomorphism of $\ov{\MM}_{2,1}(\ZZ)\setminus\ov{\WW}$
with $\ov{\UU}^{ns}_{2,1}(2)\setminus [C_0]$.
\ed

Let us consider the natural birational maps of the coarse moduli spaces
$$\ov{M}_{2,1}\dashra \ov{M}_{2,1}(\ZZ)\dashra \ov{U}^{ns}_{2,1}(2).$$
Note that all these spaces are normal (for the last two this follows from Proposition \ref{1pt-moduli-prop} and 
Corollary \ref{smooth-cor}). Note also that we only know that $\ov{M}_{2,1}(\ZZ)$ is a proper algebraic space.

Let $\ov{W}\sub\ov{M}_{2,1}$ denote the closure of $W$, and let $\De_1\sub\ov{M}_{2,1}$ be the boundary divisor,
whose generic point corresponds to the union of two elliptic curves.

\begin{prop}\label{bir-contr-prop} 
The natural birational morphism $f:\ov{M}_{2,1}\dashra \ov{M}_{2,1}(\ZZ)$ 
(resp., $g:\ov{M}_{2,1}\dashra \ov{U}^{ns}_{2,1}(2))$ is a birational contraction with the exceptional divisor $\De_1$ 
(resp., exceptional divisors $\De_1$ and $\ov{W}$).
\end{prop}

\Pf . Recall that to check that $f$ (resp., $g$) is a birational contraction we need to check
that the exceptional locus $\Exc(f^{-1})$ (resp., $\Exc(g^{-1})$)
has codimension $\ge 2$. But this immediately follows from Lemma \ref{g2-non-nodal-lem}.
Next, the restriction of $f$ to the complement of $\De_1$ induces an isomorphism with the open subset in
$\ov{M}_{2,1}(\ZZ)$ consisting of $(C,p)$ with $C$ smooth or nodal, so we have an inclusion $\Exc(f)\sub\De_1$. On the other hand, the generic point of $\De_1$
corresponds to the union of elliptic curves $E_1\cup E_2$, with the marked point on $E_1$. Under the map $f$ this curve
gets replaced by the cuspidal curve $\ov{E}_1$, so that we have a contraction $E_1\cup E_2\to \ov{E}_1$
sending the elliptic tail $E_2$ to the cusp. Since this map forgets the $j$-invariant of $E_2$, this means that $\De_1$ gets
contracted by $f$. Now the fact that $\Exc(g)=\De_1\cup\ov{W}$ follows from Theorem \ref{g2-thm}.
\ed

\begin{rem}\label{bir-contr-rem}
Let $\ov{\VV}^{ns}_{g,1}(g)\sub \ov{\UU}^{ns}_{g,1}(g)$ be the irreducible component consisting of smoothable curves.
Theorem A implies that the natural birational map
$$\ov{\MM}_{g,1}\dashra \ov{\VV}^{ns}_{g,1}(g)$$
contracts $\ov{\WW}$ to a point. Passing to the normalizations of the coarse moduli spaces we get the birational
map $\phi:\ov{M}_{g,1}\dashra X$, where $X$ is a normal projective variety, contracting $\ov{W}$ to a point.
It seems plausible that $\phi$ is a birational contraction (which would imply that $\ov{W}$ is an extremal divisor). 
To check this we would need to prove that $\Exc(\phi^{-1})$
has codimension $\ge 2$. In other words, we would need to check that the locus in $\ov{\VV}^{ns}_{g,1}(g)$,
consisting of unstable (i.e., non-nodal) curves, has codimension $\ge 2$. In the case $g=2$ we have shown this in 
Lemma \ref{g2-non-nodal-lem}. 
Note that the fact that the class of $\ov{W}$ generates an extremal ray in 
$\ov{NE}^1(\ov{M}_{g,1})$
is known for $g\le 3$ and $g=5$, by the works \cite{Rulla}, \cite{Jensen1} and \cite{Jensen2}.
\end{rem}

\end{document}